\documentclass{article}
\usepackage{ijcai20}

\usepackage{times}
\usepackage{soul}
\usepackage{url}
\usepackage[hidelinks]{hyperref}
\usepackage[utf8]{inputenc}
\usepackage[small]{caption}
\usepackage{graphicx}
\usepackage{amsmath}
\usepackage{amsthm}
\usepackage{booktabs}
\usepackage{algorithm}
\usepackage{algorithmicx}
\usepackage{amssymb,amsmath,amsfonts,graphicx,algpseudocode} 
\def\X{{\mathbb X}}
\newcommand{\lpr}{{\underline{P}}}
\newcommand{\upr}{{\overline{P}}}

\title{Uncertainty measures: The big picture}

\author{
    Fabio Cuzzolin
    \affiliations
    Visual Artificial Intelligence Laboratory, Oxford Brookes University, Oxford, UK
    \emails
    fabio.cuzzolin@brookes.ac.uk
}

\begin{document}

\maketitle

\begin{abstract}
Probability theory is far from being the most general mathematical theory of uncertainty. A number of arguments point at its inability to describe second-order (`Knightian') uncertainty. In response, a wide array of theories of uncertainty have been proposed, many of them generalisations of classical probability. As we show here, such frameworks can be organised into clusters sharing a common rationale, exhibit complex links, and are characterised by different levels of generality. Our goal is a critical appraisal of the current landscape in uncertainty theory.
\end{abstract}

\section{Introduction}

\emph{Uncertainty} is of paramount importance in artificial intelligence, applied science, and many other areas of human endevour \cite{cuzzolin2021springer}. Uncertainty can be understood as lack of information about an issue of interest for a certain agent. 
Sometimes less perceived by scientists themselves, a iatus exists between two fundamentally distinct forms of uncertainty. Somewhat `\emph{predictable}' variations (e.g., the outcomes of a carnival wheel) are typically encoded by probability distributions. \emph{Unpredictable} variations, on the other hand, reflect a more fundamental (sometimes called `Knightian') uncertainty about the laws which themselves govern the outcome (e.g., what happens if part of the wheel is cloaked?).

The mainstream mathematical theory of (first order) uncertainty is measure-theoretical probability, mainly due to Russian mathematician Andrey Kolmogorov \cite{Kolmogorov33}.
A \emph{probability measure} over a $\sigma$-algebra $\mathcal{F} = \mathcal{F}(\Theta) \subset 2^{\Theta} \doteq \{A \subseteq \Theta\}$, associated with a sample space $\Theta$, is a function $P:\mathcal{F} \rightarrow [0,1]$ such that: 
\begin{enumerate}
\item
$P(\emptyset)=0$; 
\item
$P(\Theta)=1$; 
\item
if $A\cap B = \emptyset,\; A,B\in \mathcal{F}$ then $P(A\cup B) = P(A) + P(B)$ (\emph{additivity}).
\end{enumerate}
A sample space $\Theta$ together with a $\sigma$-algebra $\mathcal{F}$ of its subsets and a probability measure $P$ on $\mathcal{F}$ forms a \emph{probability space}: $(\Theta,\mathcal{F},P)$.
A (measurable) function $X$ from a sample space $\Theta$ (endowed with a probability space) to a measurable space $E$ (usually the real line $\mathbb{R}$) is termed \emph{random variable} (RV).
Bayesian reasoning 
or frequentist inference is then applied to infer additive probabilities from the available data.

A long series of students, however, have argued that Kolmogorov's measure-theoretical probability theory is incapable of handling higher-order uncertainty, associated with widespread phenomena such as set-valued observations (which arise from missing or incomplete data), the representation of ignorance, 
inherently propositional data, and so on.
In response, various mathematical theories of uncertainty have been proposed \cite{halpern03book,maass06philosophical}. 

A theory of uncertainty is typically composed of two elements: (i) a mathematical object encoding an uncertain state of the world, and (ii) an operator which allows us to reason with uncertain states (e.g., the Choquet integral for monotone capacities, or Dempster's rule for belief functions).
The field has been (partially) surveyed by various authors \cite{Klir2004git,smets98which}, who focused, for instance, on the difference between imprecision and uncertainty, or the applicability of various models of uncertainty. 
A number of unification attempts have also been made \cite{2008arXiv0808.2747D},
most notably Walley's theory of imprecise probability \cite{walley91book,walley00towards}, Klir's generalised information theory \cite{Klir2004git}, and Zadeh's generalized theory of uncertainty \cite{Zadeh200615}. A fascinating effort to identify a unifying framework in the language of modal logic can be found in \cite{resconi93integration}.

\begin{table*} \label{tab:timeline}
\begin{center}
\caption{Imprecise-probabilistic theories: a timeline} \label{tab:timeline}
\begin{tabular}{|p{2.7cm}||p{2.6cm}|p{5cm}|p{0.8cm}|}
\hline
\footnotesize
{Approach} & \footnotesize  {Proposer(s)} & \footnotesize {Seminal paper} & \footnotesize {Year}
\\
\hline
\footnotesize
{Interval probabilities} & \footnotesize John Maynard Keynes & \footnotesize {A treatise on probability} & \footnotesize 1921
\\
\footnotesize {Subjective probability} & \footnotesize Bruno de Finetti & \footnotesize {Sul significato soggettivo della probabilit\`a} & \footnotesize 1931
\\
\footnotesize
{Theory of previsions} & 
\footnotesize
Bruno de Finetti & 
\footnotesize
{La pr\'evision: ses lois logiques, ses sources subjectives}  & 
\footnotesize
1937
\\
\footnotesize
{Theory of capacities} & 
\footnotesize
Gustave Choquet & 
\footnotesize
{Theory of capacities} & 
\footnotesize
1953
\\
\footnotesize
{Fuzzy theory} & 
\footnotesize
Lotfi Zadeh, Dieter Klaua & 
\footnotesize
{Fuzzy sets} & 
\footnotesize
1965
\\
\footnotesize
{Theory of evidence} & 
\footnotesize
Arthur Dempster, Glenn Shafer & 
\footnotesize
{Upper and lower probabilities induced by a multivalued mapping; A mathematical theory of evidence} & 
\footnotesize
1967, 1976
\\
\footnotesize
{Fuzzy measures} & 
\footnotesize
Michio Sugeno & 
\footnotesize
{Theory of fuzzy integrals and its applications} & 
\footnotesize
1974
\\
\footnotesize
{Credal sets} & 
\footnotesize
Isaac Levi & 
\footnotesize
{The enterprise of knowledge} & 
\footnotesize
1980
\\
\footnotesize
{Possibility theory} &  
\footnotesize
Didier Dubois, Henri Prade & 
\footnotesize
{Th\'eorie des possibilit\'es} & 
\footnotesize
1985
\\
\footnotesize
{Imprecise probability} & 
\footnotesize
Peter Walley & 
\footnotesize
{Statistical reasoning with imprecise probabilities} & 
\footnotesize
1991
\\
\footnotesize
{Game-theoretical probability} & 
\footnotesize
Glenn Shafer, Vladimir Vovk & 
\footnotesize
{Probability and finance: It's only a game!} & 
\footnotesize
2001
\\
\hline
\end{tabular}
\end{center}
\end{table*}

\subsection{Structure and contributions of the paper}

In this paper we will first recall the attempts to define the notion of uncertainty, and summarise the rationale behind the need to move beyond classical probability theory (Section \ref{sec:uncertainty}).
Subsequently, in the core of the paper we will review the basic concepts of each approach to uncertainty theory.

Although characterised by various levels of generality, uncertainty theories cannot be simply arranged into a `linear' hierarchy but, we argue here, should rather be
clustered according to their common rationale. 
A first group of approaches 
(Section \ref{sec:robust-probability}) 
attempt to `robustify' probability theory. Most of them are encompassed by `imprecise probability', a behavioural approach (Section \ref{sec:behavioural}) inspired by de Finetti's pioneering work. More recently, a game-theoretical approach by Vovk and Shafer has received considerable attention.
A number of frameworks, including possibility theory or rough sets and more recently Zadeh's generalised theory of uncertainty 
(Section \ref{sec:generalising-set-theory}) follow from different generalisations of set theory, upon which probability is built.
On a more general plane lie theories which generalise measure theory, the very foundation of mathematical probability (Section \ref{sec:generalised-measure}), including the theory of monotone capacities and Liu's uncertainty theory.
Some proposals (Section \ref{sec:generalised-bayesian}), such as Harper's Popperian approach and Groen's extension aim at generalising Bayesian reasoning in terms of either measures or inference mechanisms.
At the intersection of all these like the theory of random sets (or belief functions, in the finite case), which generalise the notion of set, are monotone capacities and generalise Bayesian reasoning (Section \ref{sec:set-valued}). Methods replacing events with functions in a functional space are
discussed in Section \ref{sec:others-hop}.

Due to lack of space, we neglect relatively minor proposals such as 
Cohen's {endorsements}
\cite{Cohen:1989:THR:107368.107431,Cohen:1983:FHR:1623373.1623456},
Laskey's {assumptions} \cite{blackmondlaskey89assumptions},
Shastri's evidential reasoning in semantic networks \cite{Shastri:1985:ERS:912213} or
evidential confirmation theory \cite{DBLP:journals/corr/abs-1304-3439}. 
We also entirely neglect (imprecise-)probabilistic logic \cite{ruspini86logical,Saffiotti92,Josang2001,haenni05isipta,Fagin:1990:LRP:83884.83888,wilson93-default,harmanec94modal}, which is also connected to uncertainty theory, especially in its modal logic form, but would deserve a separate survey paper. Another relevant paradigm not covered here is info-gap theory \cite{BenHaim2006}.

A critical appraisal of the uncertainty theory landscape concludes the document.


\section{Uncertainty} \label{sec:uncertainty}

\subsection{Notion of uncertainty}

\emph{Uncertainty} is of paramount importance in artificial intelligence, applied science, and many other areas of human endeavour. Whilst each and every one of us possesses some intuitive grasp of what uncertainty is, providing a formal definition can prove elusive. Uncertainty can be understood as a lack of information about an issue of interest for a certain agent (e.g., a human decision maker or a machine), a condition of limited knowledge in which it is impossible to exactly describe the state of the world or its future evolution.

According to Dennis Lindley \cite{Lindley2006}: 

``\emph{
There are some things that you know to be true, and others that you know to be false; yet, despite this extensive knowledge that you have, there remain many things whose truth or falsity is not known to you. We say that you are uncertain
about them. You are uncertain, to varying degrees, about everything in the
future; much of the past is hidden from you; and there is a lot of the present
about which you do not have full information. Uncertainty is everywhere and you
cannot escape from it}
''.

What is sometimes less clear to scientists themselves is the existence of a hiatus between two fundamentally distinct forms of uncertainty. The first level consists of somewhat `\emph{predictable}' variations, which are typically encoded as probability distributions. For instance, if a person plays a fair roulette wheel they will not, by any means, know the outcome in advance, but they will nevertheless be able to predict the frequency with which each outcome manifests itself (1/36), at least in the long run. The second level is about `\emph{unpredictable}' variations, which reflect a more fundamental uncertainty about the laws themselves which govern the outcome. Continuing with our example, suppose that the player is presented with ten different doors, which lead to rooms each containing a roulette wheel modelled by a different probability distribution. They will then be uncertain about the very game they are supposed to play. How will this affect their betting behaviour, for instance?

Lack of knowledge of the second kind is often called \emph{Knightian} uncertainty \cite{knight2012risk,Hoffman94}, from the Chicago economist Frank Knight. He would famously distinguish `risk' from `uncertainty':

``\emph{Uncertainty must be taken in a sense radically distinct from the familiar notion of
risk, from which it has never been properly separated \ldots The essential fact is that `risk' means in some cases a quantity susceptible of measurement, while at other times it is something distinctly not of this character; and there are far-reaching and crucial differences in the bearings of the phenomena depending on which of the two is really present and operating \ldots It will appear that a
measurable uncertainty, or `risk' proper, as we shall use the term, is so far different from an unmeasurable one that it is not in effect an uncertainty at all.}''

\noindent In Knight's terms, `risk' is what people normally call \emph{probability} or \emph{chance}, while the term `uncertainty' is reserved for second-order uncertainty. The latter has a measurable consequence on human behaviour: people are demonstrably averse to unpredictable
variations (as highlighted by \emph{Ellsberg's paradox} \cite{ellsberg1961risk}).

This difference between predictable and unpredictable variation is one of the fundamental issues in the philosophy of probability, and is sometimes referred to as the distinction between \emph{common cause} and \emph{special cause} \cite{snee1990statistical}. Different interpretations of probability treat these two aspects of uncertainty in different ways, as debated by economists such as John Maynard Keynes \cite{keynes1922treatise} and G. L. S. Shackle.

\subsection{Probability}

Measure-theoretical probability, due to the Russian mathematician Andrey Kolmogorov \cite{Kolmogorov33}, is the mainstream mathematical theory of (first-order) uncertainty. In Kolmogorov's mathematical approach  probability is simply an application of measure theory, and uncertainty is modelled using additive measures.

A number of authors, however, have argued that measure-theoretical probability theory is not quite up to the task when it comes to encoding second-order uncertainty. In particular, as we discuss in the Introduction, additive probability measures cannot properly model missing data or data that comes in the form of \emph{sets}.
Probability theory's frequentist interpretation is utterly incapable of modelling `pure' data (without `designing' the experiment which generates it). In a way, it cannot even properly model continuous data (owing to the fact that, under measure-theoretical probability, every point of a continuous domain has zero probability), and has to resort to the `tail event' contraption to assess its own hypotheses. Scarce data can only be effectively modelled asymptotically.

Bayesian reasoning is also plagued by many serious limitations: (i) it just cannot model ignorance (absence of data); (ii) it cannot model pure data (without artificially introducing a prior, even when there is no justification for doing so); (iii) it cannot model `uncertain' data, i.e., information not in the form of propositions of the kind `$A$ is true'; and (iv) again, it is able to model scarce data only asymptotically, thanks to the Bernstein--von Mises theorem \cite{cramer2016mathematical}.

\subsection{Beyond probability}

Similar considerations have led a number of scientists to recognise the need for a coherent mathematical theory of uncertainty able to properly tackle all these issues.
Both alternatives to and extensions of classical probability theory have been proposed, starting from de Finetti's pioneering work on subjective probability \cite{DeFinetti74}. 
Formalisms include possibility-fuzzy set theory \cite{Zadeh78,Dubois90}, probability intervals
\cite{halpern03book}, credal sets \cite{levi80book,kyburg87bayesian}, monotone capacities \cite{wang97choquet}, random sets \cite{Nguyen78} and imprecise probability theory \cite{walley91book}. New original foundations of subjective probability in behavioural terms \cite{walley00towards} or by means of game theory \cite{shafer01book} have been put forward.

Table \ref{tab:timeline} presents a sketchy timeline of the various existing approaches to the mathematics of uncertainty.
Sometimes collectively referred to as \emph{imprecise probabilities} (as most of them comprise classical probabilities as a special case), these theories in fact form an entire hierarchy of encapsulated formalisms.

In what follows we will attempt to classify those formalisms into clusters characterised by a common theme.


\section{Probability robustified} \label{sec:robust-probability}

\subsection{Lower and upper probabilities} \label{sec:lower-upper}

A \emph{lower probability} (LP) \cite{FINE1988389} $\lpr$ is a function from a sigma-algebra 
to the unit interval $[0,1]$ such that:
\[
\underline{P}(A \cup B) \geq \underline{P}(A) + \underline{P}(B) \quad \forall A \cap B = \emptyset
\]
(super-additivity). $\overline{P}$ is an \emph{upper probability} (UP) whenever 
\[
\overline{P}(A \cup B) \leq \overline{P}(A) + \overline{P}(B), 
\quad
A \cap B = \emptyset.
\]
With a lower probability $\lpr$ is associated a dual UP $\overline{P}(A)=1-\underline{P}(A^c)$, $A^c$ being the complement of $A \subseteq \Theta$.

\subsection{Probability envelopes or credal sets}

Each $\lpr$ is also associated with a closed convex set of probability measures 
\[
\mathcal{P}(\lpr)= \{ P : P(A) \geq \lpr(A), \forall A \subseteq \Theta \}, 
\]
termed a \emph{credal set} \cite{levi80book}.
Not all credal sets, however, can be described by merely focusing on events \cite{walley91book}.

A lower probability $\lpr$ on $\Theta$ is called `consistent' (`avoids sure loss' in Walley's terminology) if $\mathcal{P}(\lpr)\neq \emptyset$, i.e.:
\[
\sup_{\theta \in\Theta} \sum_{i=1}^n 1_{A_i}(\theta) \geq \sum_{i=1}^n \lpr(A_i),
\]
whenever $n \in \mathbb{N}^+$, $A_1,...,A_n \in \mathcal{F}$ are events, 
and $1_{A_i}$ is the indicator function of $A_i$ in $\Theta$.
\\
A lower probability $\lpr$ is `tight' (`coherent' for Walley) if $\inf_{P \in \mathcal{P}(\lpr)} P(A) = \lpr(A)$, i.e.:
\[
\sup_{\theta\in\Theta} [ \sum_{i=1}^n 1_{A_i}(\theta) - m 1_{A_0}(\theta) ] \geq \sum_{i=1}^n \lpr(A_i) - m \cdot \lpr(A_0) 
\]
whenever $n,m \in \mathbb{N}^+$ and $A_0,A_1,...,A_n \in \mathcal{F}$.
Consistency means that the lower bound constraints $\lpr(A)$ can indeed be satisfied by some probability measure, while tightness indicates that $\lpr$ is the lower envelope on subsets of $\mathcal{P}(\lpr)$. Any coherent lower probability is \emph{monotone}, i.e., $\lpr(A)\leq \lpr(B)$ for $A \subset B$.

As shown by \cite{kyburg87bayesian}, one can operate with credal sets in a robust statistical fashion, by simply applying Bayes' rule to their (finite) set of vertices.

\subsection{Probability intervals} \label{sec:soa-others-intervals} \label{sec:others-intervals}

One can instead provide lower and upper constraints on the probabilities of \emph{elements} $x \in \Theta$.
A \emph{set of probability intervals} \cite{Kyburg98interval-valuedprobabilities,tessem92interval,decampos94} is a system of constraints 
\[
\mathcal{P}(l,u) \doteq \{p: l(x) \leq p(x) \leq u(x), \forall x \in \Theta \}
\]
on the values of a probability distribution $p:\Theta \rightarrow [0,1]$ on a finite domain $\Theta$.
Such a set of constraints determines a credal set, of a sub-class of those generated by lower and upper probabilities.
De Campos et al. studied the specific constraints such intervals need satisfy in order to be consistent and tight. Their vertices can be computed as in \cite{decampos94}, p. 174.

In particular, a set of probability intervals is \emph{feasible} if and only if for each
$x \in \Theta$ and every value $v(x) \in [l(u), u(x)]$ there exists a probability distribution function $p:\Theta \rightarrow [0,1]$ for which $p(x) = v(x)$.
One can then obtain the lower and upper probabilities on any subset $A \subseteq \Theta$ via: 
\[
\begin{array}{l}
\displaystyle
\lpr(A) = \max \bigg \{ \sum_{x \in A} l(x), 1- \sum_{x \not\in A} u(x) \bigg \},  
\\
\displaystyle
\lpr(A) = \min \bigg \{ \sum_{x \in A} u(x), 1- \sum_{x \not\in A} l(x) \bigg \}.
\end{array}
\]
Lower and upper probabilities associated with a feasible set of probability intervals are 2-monotone capacities (\cite{decampos94}, Prop 5; see Section \ref{sec:others-capacities}).
Belief functions (Section \ref{sec:belief-functions}) also form a special class of interval probabilites. 

Combination, marginalisation and conditioning operators for probability intervals can be defined, typically using operators acting on lower and upper probabilities. In \cite{decampos94}, \emph{credal conditioning} \cite{fagin91new} was suggested: 
\[
\lpr(A|B) = \frac{\lpr(A \cap B)}{\lpr(A \cap B) + \upr(A^c \cap B)} 
\quad
\forall A, B \subseteq \Theta.
\]
A generalised Bayesian inference framework based on interval probabilities was also proposed in \cite{pan97bayesian}.

\subsection{Probability boxes} \label{sec:soa-others-pboxes}  \label{sec:others-pboxes}

\emph{Probability boxes} \cite{Ferson03pboxes,ALVAREZ2006241,joslyn04approximate} arise in engineering and reliability whenever the available information is insufficient to identify a sought joint probability density function (PDF).
Let $P$ be a probability measure on the real line $\mathbb{R}$. Its
\emph{cumulative distribution function} (CDF) is a non-decreasing mapping $F_P$ from $\mathbb{R}$ to $[0,1]$, such that for any $\alpha \in \mathbb{R}$, $F_P(\alpha) = P((-\infty, \alpha])$. 

A \emph{probability box} or \emph{p-box} \cite{Ferson03pboxes} $\langle \underline{F},\overline{F} \rangle$ is a class of cumulative distribution functions
delimited by two lower and upper bounds $\underline{F}$ and $\overline{F}$:
\[
\langle \underline{F},\overline{F} \rangle = \{ F \; \text{CDF} |  \underline{F} \leq F \leq \overline{F} \}.
\]
P-boxes can be convolved by first discretising them into belief functions (cfr. Section \ref{sec:belief-functions}) and applying a \emph{parametric copula} \cite{Yager2013} to the latter. 
A copula $C:[0,1]^d \rightarrow [0,1]$ is a multivariate probability distribution for which the marginal probability distribution of each variable is uniform. Copulas are used to describe the dependence between RVs.
The resulting BF is finally transformed back into a p-box.

Generalised p-boxes \cite{2008arXiv0808.2747D} have also been proposed.

\section{Behavioural probability} \label{sec:behavioural}

In behavioural probability, pioneered by \cite{DeFinetti74,deFinetti1980}, probability has a behavioural rationale which derives from equalling `belief' to `inclination to act'. An agent believes in a certain outcome to the extent they are willing to accept a bet on that outcome.

\subsection{De Finetti's previsions} 

A person who wants to summarise his degree of belief in a random event $A \subset \Theta$ by a number $p$ is supposed to accept any bet on $A$ with gain $c(p - 1_A)$, where $1_A$ denotes the indicator function of $A$ and $c$ is any real number chosen by an opponent.
Since $c$ may be positive or negative, there is no advantage for the person in question in choosing a value $p$ such that $cp$ and $c (p-1)$ are both strictly positive for some $c$. Hence, $p$ is an admissible evaluation of the probability of $A$ if it meets a \emph{principle of coherence}: $p$ has to ensure that $\nexists c \in \mathbb{R}$ so that the realisations of $c(p -1_A)$ are all strictly positive (or strictly negative). 
Coherence means that the agent is fully aware of the consequences of its betting rates.

\subsubsection{Probabilities}

The concept of coherence is extended to a class $\mathcal{F}$ of events as follows. The real-valued function $P$ on $\mathcal{F} \subset 2^\Theta$ is said to be a
\emph{probability} if, for any finite subclass $\{A_1,...,A_n\}$ of $\mathcal{F}$ and any choice of $(c_1,...,c_n)$ in $\mathbb{R}^n$, $n=1,2,...$ the gain 
\[
G = \sum_{k=1}^n (P(A_k) - 1_{A_k}) 
\]
is such that $\inf G \leq 0 \leq \sup G$, where $\inf$ and $\sup$ are taken over all constituents of $\{A_1,..., A_n\}$.

\subsubsection{Gambles and previsions}

Now, let a \emph{gamble} $X \in \mathcal{L}(\Theta)$ \cite{DeCooman2003agentle} 
be a bounded real-valued function on $\Theta$: $X:\Theta \rightarrow \mathbb{R}$, $\theta \mapsto X(\theta)$ yielding different utilities for different outcomes $\theta \in \Theta$. 

$P$ is termed a
\emph{prevision} if, for every finite subclass $\{X_1,...,X_n\}$ of $\mathcal{L}(\Theta)$ and for every $n$-tuple $(c_1,...,c_n)$ of real numbers, one has
\[
\inf \sum_{k=1}^n c_k (P(X_k) - X_k) \leq 0 \leq \sup \sum_{k=1}^n c_k (P(X_k) - X_k),
\]
i.e., there is no (finite) betting system which makes uniformly strictly negative the gain of the agent adopting $P$. The probability of an event $A$ coincides with the prevision of $1_A$ and the theory of probability is included in that of previsions.

A prevision $P$ is \emph{coherent} iff
\begin{enumerate}
\item
$P(\lambda X + \mu Y) = \lambda P(X) + \mu P(Y)$; 
\item
if $X>0$ then $P(X)\geq 0$; 
\item
$P(\Theta) = 1$.
\end{enumerate}

\subsection{Imprecise probability} \label{sec:others-imprecise} \label{sec:walley}

Inspired by de Finetti's work, \emph{imprecise probability} \cite{walley00towards,walley91book} (IP, see \cite{miranda2008survey} for a survey) aims at unifying all approaches to mathematical uncertainty in a single coherent setting, with a generality comparable to the theory of credal sets \cite{Cozman2000191}. 

\subsubsection{Sets of desirable gambles}

There, an agent's \emph{set of desirable gambles} $\mathcal{D} \subseteq \mathcal{L}(\Theta)$ is used as a model of their uncertainty about the problem.
$\mathcal{D}$ is \emph{coherent} iff it is a convex cone, i.e.,
\begin{enumerate}
\item
0 (the constant gamble $X(\theta) = 0$ $\forall\theta$) $\not\in \mathcal{D}$;
\item
if $X>0$ (i.e., $X(\theta) > 0$ for all $\theta$) then $X \in \mathcal{D}$;
\item
if $X,Y \in\mathcal{D}$, then $X + Y \in \mathcal{D}$;
\item
if $X \in\mathcal{D}$ and $\lambda >0$ then $\lambda X \in \mathcal{D}$.
\end{enumerate}
As a consequence, if $X \in\mathcal{D}$ and $Y > X$ then $Y \in\mathcal{D}$.

\subsubsection{Lower and upper previsions}

Suppose the agent buys a gamble $X$ for a price $\alpha$. This yields a new gamble $X - \alpha$.

The \emph{lower prevision} $\underline{P}(X)$ of a gamble $X$,
\[
\underline{P}(X) \doteq \sup \{\alpha : X - \alpha \in \mathcal{D} \},
\]
is the supremum acceptable price for buying $X$.
Selling a gamble $X$ for a price $\alpha$ yields a new gamble $\alpha - X$.

The \emph{upper prevision} 
\[
\overline{P}(X) \doteq \inf \{\alpha : \alpha - X \in \mathcal{D} \}
\]
is the supremum acceptable price for selling $X$.
By definition, $\overline{P}(X) =  - \underline{P}(-X)$. When lower and upper previsions coincide,
$\overline{P}(X) =  \underline{P}(X) = P(X)$ is called the (precise) \emph{prevision} of $X$ in de Finetti's sense \cite{DeFinetti74}.

Rational rules of behaviour such as `avoiding sure loss' and coherence can be applied to lower previsions too, by replacing indicator functions with general gambles.
One consequence of avoiding sure loss is that $\underline{P}(A) \leq \overline{P}(A)$. From coherence it follows that lower previsions are subadditive.

\subsubsection{Natural and marginal extension}

\emph{Natural extension} is used to extend a coherent lower prevision defined on a collection of gambles to a lower prevision on all gambles, assumed coherent and conservative. 

The natural extension of a set of gambles $\mathcal{D}$ is the smallest coherent set of desirable gambles that includes $\mathcal{D}$. 
Let $\lpr$ be a lower probability on $\Theta$ that avoids sure loss, and let $\mathcal{L}$ be the set of all bounded functions on $\Theta$. The functional $\underline{E}$ defined on $\mathcal{L}$ as $\underline{E}(f) =$
\[
\begin{array}{l}
\displaystyle
 \sup 
\bigg \{ 
\sum_{i=1}^n \lambda_i \lpr(A_i) + c \; | \; f \geq \sum_{i=1}^n \lambda_i 1_{A_i} + c, n \geq 0, 
\\
\displaystyle
A_i \subseteq \Theta, \lambda_i \geq 0, c \in \{-\infty, +\infty\} 
\bigg \}
\end{array}
\] 
for all $f \in \mathcal{L}$ is called the \emph{natural extension} of $\lpr$.
A similar definition can be given for lower previsions.
When $\lpr$ is a classical (`precise') probability, the natural extension agrees with the expectation. Also, $\underline{E}(1_A) = \lpr(A)$ for all $A$ iff $\lpr$ is coherent.

A \emph{marginal extension} operator \cite{miranda2008survey} can be introduced for the aggregation of conditional lower previsions.

\subsection{Game-theoretical probability} \label{sec:others-vovk} \label{sec:vovk}

In \emph{game-theoretic probability} \cite{shafer01book} 
probabilistic predictions are proven by constructing a betting strategy which allows to multiply capital indefinitely if the prediction fails, in a game-theoretical setting. Probabilistic theories are tested by betting against its predictions.

Mathematically, the approach builds on the predating theory of \emph{prequential probability} \cite{dawid1999prequential}. Dawid's prequential principle states that the forecasting success of a probability distribution for a sequence of events should be evaluated using only the actual outcomes and the sequence of forecasts (conditional probabilities) to which these outcomes give rise, without reference to other aspects. In fact, the notion goes back to the work of Blaise Pascal, (`Probability is about betting'), Antoine Cournot (`Events of small probability do not happen'), and Jean Ville \cite{Shafer07game}. 

The latter, in particular, showed in 1939 that the laws of probability can be derived from this principle: \emph{You will not multiply the capital you risk by a large factor} (the so-called \emph{Cournot's principle}).
The Ville/Vovk perfect-information protocol for probability can be stated in terms of a game involving three players: Forecaster, Skeptic and Reality, as in Algorithm \ref{alg:game}.
Using Shafer and Vovk's weather forecasting example, Forecaster may be a very complex computer program that escapes precise mathematical definition as it is constantly under development, which uses information external to the model to announce every evening a probability for rain on the following day. Skeptic decides whether to bet for or against rain and how much, and Reality decides whether it rains.
As shown by Ville, the fundamental hypothesis that Skeptic cannot get rich can be tested by any strategy for betting at Forecaster’s odds.

The bottom line of Shafer and Vovk's theory if that Forecaster can beat Skeptic \cite{vovk2005defensive}.

\begin{algorithm} 
\caption{Game theoretical probability protocol \label{alg:game}}  
\begin{algorithmic}[1]
\Procedure{GameTheoreticalProtocol}{}
\State $\mathcal{K}_0 = 1$
\For{$n = 1, \cdots, N$}
\State
Forecaster announces prices for various payoffs.
\State
Skeptic decides which payoffs to buy.
\State
Reality determines the payoffs.
\State
Skeptic's capital changes as $K_n = K_{n-1} +$ net gain or loss.
\EndFor
\EndProcedure
\end{algorithmic}
\end{algorithm}

\section{Generalising set theory} \label{sec:generalising-set-theory}

\subsection{Fuzzy and possibility theory} \label{sec:others-fuzzy}

The concept of \emph{fuzzy set} was introduced in \cite{zadeh65fuzzysets} 
as an extension of the notion of set. 

\subsubsection{Possibility theory}

\emph{Possibility theory} \cite{dubois88possibility} equips fuzzy set theory to provide a graded semantics to natural language statements.

A \emph{possibility measure} on $\Theta$ is a function $\Pi: 2^\Theta \rightarrow [0,1]$ such that 
\begin{enumerate}
\item
$\Pi(\emptyset) =0$; 
\item
$\Pi(\Theta)=1$;
\item
the condition 
\[
\Pi \bigg ( \bigcup_i A_i  \bigg ) = \sup_i \Pi(A_i)
\]
is satisfied for every family of subsets $\{ A_i\in 2^\Theta \}$.
\end{enumerate}
Each possibility measure is uniquely characterised by a \emph{membership function} 
$\pi: \Theta \rightarrow [0,1]$ s.t. $\pi(x)\doteq \Pi(\{x\})$ via the formula
\[
\Pi(A) = \sup_{x\in A} \pi(x).
\] 
The dual quantity $N(A) = 1 - \Pi(A^c)$ is called \emph{necessity measure}.
A further extension called \emph{vague set} \cite{229476} imposes a pair of lower and upper bounds on the membership function of a fuzzy set.

A conditional possibility measure $\pi(.|A)$ can be defined such that $\Pi(A \cap B) = \min \{ \Pi(B|A), \Pi(A) \}$. Its least specific solution is $\Pi(B|A) = 1$ if $\Pi(A \cap B) = \Pi(A)$, $\Pi(A \cap B)$ else. Other solutions \cite{dubois97bayesian} are possible.

\subsubsection{Other fuzzy-inspired frameworks}

Other uncertainty frameworks based on fuzzy theory have been proposed, including \emph{fril-fuzzy} \cite{baldwin95-frilfuzzy}, 
\emph{granular computing} \cite{Yao00granularcomputing}, \emph{interval fuzzy reasoning} \cite{Yao97interval} and \emph{neighborhood systems} \cite{lin96Neighborhoods}.

\subsection{Rough sets} \label{sec:soa-others-rough} \label{sec:others-rough-sets} \label{sec:rough-sets}

In opposition, 
\emph{rough sets} \cite{Pawlak1982} are 
strongly linked to the idea of partition of the universe of hypotheses.

Let $\Theta$ be a finite universe, and $\mathcal{R} \subseteq \Theta \times \Theta$ be an equivalence relation  which partitions it into a family of disjoint subsets $\Theta/\mathcal{R}$, called \emph{elementary sets}.
We can then call \emph{measurable} sets $\sigma(\Theta/\mathcal{R})$ the unions of one or more elementary sets, plus the empty set $\emptyset$.
The \emph{lower approximation} $\underline{apr}(A)$ of $A \subset \Theta$ is formed by the measurable elements whose equivalence classes are contained in $A$. Its \emph{upper approximation} $\overline{apr}(A)$ is composed by the elements whose equivalence classes have nonempty intersection with $A$. 

Rough sets are useful to describe the probabilistic properties of decision rules \cite{doi:10.1080/019697298125470},
and 
are quite related to belief functions.
Given a $\sigma$-algebra $\mathcal{F}$ of subsets of $\Theta$, one can construct a rough set algebra such that $\mathcal{F} = \sigma(\Theta/\mathcal{R})$. Any probability $P$ on $\mathcal{F}$ can be extended to $2^\Theta$ using inner and outer measures as follows:
\[
\begin{array}{l}
\displaystyle
P_*(A) = \sup \{ P(X) | X \in \sigma(\Theta /\mathcal{R}), X \subseteq A \}, 
\\
\displaystyle
P^*(A) = \sup \{ P(X) | X \in \sigma(\Theta /\mathcal{R}), X \supseteq A \}.
\end{array}
\]
These are a pair of belief and plausibility functions (Section \ref{sec:belief-functions}).

\subsubsection{Generalized theory of uncertainty (GTU)} \label{sec:soa-others-gtu} \label{sec:others-gtu}

\emph{Generalised theory of uncertainty} (GTU) \cite{Zadeh200615} describes information via \emph{generalised constraints} 
of the form $GC(X) : X \text{isr} R$, where $r \in \{$ blank, probabilistic, veristic, random set, fuzzy graph, etc $\}$ is a label which determines the type of constraint, and $R$ a constraining relation of that type (e.g. a probability distribution, a random set, etc).
This allows us to operate on information such as: `Usually Robert returns from work at about 6 p.m.', 
to which generalised constraint propagation is employed as a reasoning mechanism. 

In the GTU everything is or is allowed to be a matter of degree (fuzzy).
A \emph{generalised constraint language} is defined as the set of all generalised constraints together with the rules governing syntax, semantics and generation. Examples are: ($X$ is small) is likely; $((X, Y)\; \text{isp} \; A) \wedge (X \; \text{is} \; B)$, where `isp' denotes a probabilistic constraint, `is' a possibilistic constraint, and $\wedge$ conjunction.
In GTU, inference is treated as an instance of question-answering. Given a system of natural language propositions and a query $q$, 
GTU performs generalised constraint propagation governed by deduction rules. 

\section{Generalising measure theory} \label{sec:generalised-measure}

\subsection{Capacities (fuzzy measures)} \label{sec:soa-others-fuzzy-measures} \label{sec:others-capacities}

The \emph{theory of capacities} or 
\emph{fuzzy measure theory} \cite{wang1992book,sugeno74fuzzy,grabisch2000book} 
generalises measure theory by replacing additivity by \emph{monotonicity} \cite{Choquet53,sugeno74fuzzy}. 

Many uncertainty measures are special cases of fuzzy measures, including belief functions, possibilities and probabilities \cite{LAMATA1989243}.

\subsubsection{Capacities}
Given a domain $\Theta$ and a non-empty family $\mathcal{F}$ of subsets of $\Theta$, a \emph{monotone capacity} or \emph{fuzzy measure} $\mu$ on
$(\Theta, \mathcal{F})$ is a function $\mu : \mathcal{F} \rightarrow [0,1]$ such that 
\begin{enumerate}
\item
$\mu(\emptyset)=0$;
\item
if $A\subseteq B$ then 
$
\mu(A)\leq \mu(B) 
$
(monotonicity);
\end{enumerate}
Sometimes the following additional conditions are imposed: 
\begin{enumerate}
\setcounter{enumi}{2}
\item
for any increasing sequence $A_1\subseteq A_2 \subseteq \cdots$ in $\mathcal{F}$, if $\bigcup_{i=1}^{\infty}A_i\in\mathcal{F}$ then 
\[
\lim_{i\rightarrow \infty} \mu(A_i) = \mu \bigg ( \bigcup_{i=1}^{\infty}A_i \bigg ) 
\]
(`{continuity from below}');
\item
for any decreasing sequence $A_1\supseteq A_2 \supseteq \cdots$, if
$\bigcap_{i=1}^{\infty}A_i\in\mathcal{F}$ and $\mu(A_1)<\infty$, then 
\[
\lim_{i\rightarrow \infty}\mu(A_i) = \mu \bigg ( \bigcap_{i=1}^{\infty}A_i \bigg ) 
\]
(`{continuity from above}').
\end{enumerate}
When $\Theta$ is finite the last two requirements are trivially satisfied and can be disregarded. Monotone \emph{decreasing} measures can also be defined.

\subsubsection{Choquet integral}

For any nonnegative measurable function $f$ on $(\Theta, \mathcal{F})$, the Choquet integral of $f$ on any $A \in \mathcal{F}$ is defined as
\[
\inf_A f d\mu \doteq \int_0^\infty \mu(F_\alpha \cap A) d\alpha, 
\]
where $F_\alpha = \{x \in \Theta | f(x)\geq \alpha \}$, $\alpha \in [0,\infty)$.
Both the Choquet integral of monotone capacities and the natural extension of lower probabilities are generalisations of the Lebesgue integral with respect to $\sigma$-additive measures. The Choquet integral is monotone in both arguments $f$ and $\mu$, but is not a linear functional.

\subsubsection{Order of capacities}

A capacity $\mu$ is termed of \emph{order} $k$ if 
\[
\mu \bigg ( \bigcup_{j=1}^k A_j \bigg ) 
\geq 
\sum_{\emptyset \neq K \subseteq [1,...,k]} (-1)^{|K|+1} 
\mu \bigg ( \bigcap_{j \in K} A_j \bigg ) 
\]
for all collections of $k$ subsets $A_j$, $j \in K$ of $\Theta$.
If $k'>k$ the resulting theory is less general than a theory of capacities of order $k$ (i.e., it contemplates fewer measures). A capacity is \emph{$k$-alternating} when its dual is $k$-monotone.
When a capacity is $k$-alternating (monotone) for every $k$, it is called $\infty$-alternating (monotone).

\emph{Monotone capacities of order 2} are those such that
\[
\mu(A \cup B) + \mu(A \cap B) \geq \mu(A) + \mu(B) 
\quad 
\forall A,B \subseteq \Theta. 
\]
Their conjugate measures $\mu^c \doteq 1 - \mu(A^c)$ are such that
\[
\mu^c(A \cup B) + \mu^c(A \cap B) \leq \mu^c(A) + \mu^c(B).
\]
A 2-monotone capacity is a lower probability (but not viceversa). Feasible probability intervals are 2-monotone capacities.

\subsubsection{Sugeno $\lambda$-measures}

Sugeno \emph{$\lambda$-measures} \cite{sugeno74fuzzy}, are
special monotone measures $g_\lambda$ such that
\[
g_\lambda(A \cup B) = g_\lambda(A) + g_\lambda(B) + \lambda g_\lambda(A) g_\lambda(B)
\]
for any given pair of disjoint sets $A, B \in 2^\Theta$, where $\lambda \in (-1, \infty)$ is a parameter.

Each $\lambda$-measure is uniquely determined by values $g_\lambda(\theta)$, $\theta \in \Theta$
\cite{wang1992book}:
\begin{itemize}
\item
if $\sum_\theta g_\lambda(\theta) < 1$, then $g_\lambda$ is a lower probability; 
\item
if $\sum_\theta g_\lambda(\theta) = 1$ then $g_\lambda$ is a probability measure and $\lambda = 0$; 
\item
if $\sum_\theta g_\lambda(\theta) > 1$, $g_\lambda$ is an upper probability.
\end{itemize}
Lower and upper probabilities based on $\lambda$-measures are special belief and plausibility measures \cite{Berres1988159}. 

\subsection{Comparative probability}

Probability measures can also be generalised by replacing numerical values with order relations. Fine's \emph{comparative probability} (CP) \cite{FINE197315} asserts that, given any two propositions $A$ and $B$, all that is necessary is to say is whether $P(A)$ is greater than, lower than or equal to $P(B)$. A CP relation is required: (i) to be a linear, complete and simple order; (ii) to be non-trivial ($\Theta \succ \emptyset$); (iii) to satisfy the `improbability of impossibility' condition, $A \succ \emptyset$ for all $A$; (iv) to be such that $A \cap (B \cup C) = \emptyset$ implies $B \succ C$ iff $A \cup B \succ A \cup C$.

Wong et al. later generalised this to belief structures \cite{WONG1992123}, and showed that within this system there are belief functions which `{almost agree}' with comparative relations, i.e.,
\[
A \succ B \Rightarrow Bel(A) > Bel(B) 
\quad 
\forall A,B \in 2^\Theta.
\]

\subsection{Baoding Liu's Uncertainty Theory} \label{sec:soa-others-liu} \label{sec:others-ut}

Liu's \emph{uncertainty theory} (UT) \cite{liu2004book,liu2009some} is based on the notion of \emph{uncertain measure}, a function $\mathcal{M}$ on $\mathcal{F}(\Theta)$ such that: 
\begin{enumerate}
\item
$\mathcal{M}(\Theta) = 1$ (normality); 
\item
$\mathcal{M}(A) \leq \mathcal{M}(B)$ if $A \subset B$ (monotonicity); 
\item
$\mathcal{M}(A) + \mathcal{M}(A^c) = 1$ (`self-duality'); 
\item
for every countable sequence of events $\{A_i\}$ we have that
\[
\mathcal{M}\left (\bigcup_{i=1}^\infty A_i \right ) \leq \sum_{i=1}^\infty \mathcal{M}(A_i)
\] 
(countable subadditivity).
\end{enumerate}

Clearly, uncertain measures are monotone capacities. Just as clearly, probability measures do satisfy these axioms. However, Liu claims, probability theory is not a special case of UT since probabilities do not satisfy the \emph{product axiom}:
\[
\mathcal{M}\left ( \prod_{k=1}^\infty A_k \right ) = \bigwedge_{k=1}^\infty \mathcal{M}_k (A_k)
\]
for all Cartesian products of events from individual uncertain spaces $(\Theta_k,\mathcal{F}_k,\mathcal{M}_k)$. However, the product axiom was only introduced in \cite{liu2009some} (much after Liu's introduction of uncertain theory). 
The extension of uncertain measures to any subset of a product algebra is rather cumbersome and unjustified (Equation (1.10) or Figure 1.1 in \cite{liu2004book}). Axioms are not well justified.
Nevertheless, based on such measures a generalisation of random variables can be defined (`uncertain variables'), as measurable (in the usual sense) functions from an uncertainty space $(\Theta, \mathcal{F},\mathcal{M})$ to the set of real numbers.

\section{Generalising Bayesian reasoning} \label{sec:generalised-bayesian}

\subsection{Harper's Popperian approach}

\cite{10.2307/20115029} proposed 
to use Popper's \emph{propensity} treatment of probability 
to extend Bayesian reasoning so that revision of previously accepted evidence is allowed for.

A function $P: \mathcal{F} \times \mathcal{F} \rightarrow [0,1]$ assigning values to pairs of 
events, endowed with a binary operation $AB$ and a unary one $\bar{A}$, is a \emph{Popper probability function} if:
\begin{enumerate}
\item $0 \leq P(B|A) \leq P(A|A) =1$;
\item if $P(A|B) = 1 = P(B|A)$ then $P(C|A) = P(C|B)$;
\item if $P(C|A) \neq 1$ then $P(\bar{B}|A) = 1 - P(B|A)$;
\item $P(AB|C) = P(A|C) \cdot P(B|AC)$;
\item $P(AB|C) \leq P(B|C)$.
\end{enumerate}

In Popperian probability, conditional probability is primitive, whereas absolute one is derived as 
\[
P(A) \doteq P(A|T), 
\quad
T = \overline{A\bar{A}}.
\]
Note that on traditional $\sigma$-algebras $T = (A \cap A^c)^c = \Theta$. When $P(A)$, as defined above, is greater than 0 the Popper conditional probability is the ratio of the absolute probabilities: $P(B|A) = P(AB)/P(B)$, i.e, it is a classical conditional probability. 
In Popper's theory, however, $P(B|A)$ exists even when $P(A)=0$.
As a result, previously accepted evidence can be revised to condition on events which currently have 0 probability.

Harper's results included an epistemic semantics for the theory of \emph{counterfactual conditionals} \cite{Lewis1971,Lewis1973}.

\subsection{Groen's extension of Bayesian theory}

\cite{GROEN200549} proposed an extension of Bayesian theory in which observations are used to rule out possible valuations of the variables.
The extension is different from probabilistic approaches such as Jeffrey's rule, in which certainty in a single proposition $A$ is replaced by a probability on a disjoint partition of the universe, and Cheeseman's rule of distributed meaning \cite{Cheeseman86}, while non-probabilistic analogues are found in evidence and possibility theory. 

\subsubsection{Interpretations}

An \emph{interpretation} $I$ of an observation $O$ is defined as the union of values of a variable of interest $X$ (e.g., an object's luminosity) that are not contradicted by $O$ (e.g., the object is `very dark'). Uncertainty regarding the interpretation of an observation can be expressed as either a PDF on the space $\mathcal{I}$ of possible interpretations, $\pi(I|H)$, $I \in \mathcal{I}$,
or by introducing an \emph{interpretation function}
\[
\rho(x) \doteq Pr(x \in I|H), 
\quad
0 \leq \rho(x) \leq 1, 
\]
where $x \in I$ denotes the event that value $x$ is not contradicted by the observation,
$\rho(x)$ is the probability that this is true, and $H$ our prior knowledge.
\\
The two definitions are related via:
\[
\rho(x) = 
\sum_{I \in \mathcal{I}} Pr(x \in I|I, H) \pi(I|H) = \sum_{I \in \mathcal{I}: x \in I} \pi(I|H).
\]

Prior to making an observation, our initial state of uncertainty about the pair `observable variable' $x \in \X$, `variable of interest' $\theta \in \Theta$, can be represented by the probability distribution:
$\pi(x,\theta|H) = \pi(x| \theta, H) \pi(\theta|H)$ defined over $\X \times \Theta$.
Full Bayesian inference can then be modified by introducing the above interpretation function into Bayes' rule, yielding the posterior joint: 
\[
\pi(x,\theta|H,O) = \frac{\rho(x) \pi(x|\theta,H) \pi(\theta|H)}{\int_{x,\theta} \rho(x) \pi(x|\theta,H) \pi(\theta|H) dx d\theta}.
\]
The rationale is that observations provide no basis for a preference among representations not contradicted by them. Thus, the relative likelihoods of the remaining representations should not be affected by it.

\section{Set-valued probability}  \label{sec:set-valued}

\subsection{Belief functions} \label{sec:belief-functions}

Let us denote by $\Omega$ and $\Theta$ the sets of outcomes of two different but related problems $Q_1$ and $Q_2$, respectively. Given a probability measure $P$ on $\Omega$, we want to derive a `degree of belief' $Bel(A)$ that $A\subset \Theta$ contains the correct response to $Q_2$. If we call $\Gamma(\omega)$ the subset of outcomes of $Q_2$ compatible with $\omega\in\Omega$, $\omega$ tells us that the answer to $Q_2$ is in $A$ whenever $\Gamma(\omega) \subset A$. The \emph{degree of belief} $Bel(A)$ of an event $A\subset\Theta$ is then the total probability (in $\Omega$) of all the outcomes $\omega$ of $Q_1$ that satisfy the above condition, namely \cite{Dempster67}:
\[
Bel(A) = P(\{ \omega | \Gamma(\omega) \subset A \}) = \sum_{\omega\in\Omega
: \Gamma(\omega)\subset A} P(\{\omega\}).
\]
The map $\Gamma : \Omega \rightarrow 2^{\Theta} = \{A \subseteq \Theta\}$ 
is called a \emph{multivalued mapping} 
from $\Omega$ to $\Theta$. Such a mapping, together with a probability measure $P$ on $\Omega$, induces a \emph{belief function} on $2^\Theta$.

\subsubsection{Belief and plausibility measures}

A \emph{basic probability assignment} (BPA) \cite{Shafer76} 
is a set function \cite{denneberg99interaction,dubois86logical} 
$m : 2^\Theta\rightarrow[0,1]$ 
s.t.: 
\begin{enumerate}
\item
$m(\emptyset)=0$; 
\item
$\sum_{A\subset\Theta} m(A)=1$.
\end{enumerate}
The `mass' $m(A)$ 
assigned to $A$ is in fact the probability $P(\{\omega \in \Omega : \Gamma(\omega) = A\})$.
Subsets of $\Theta$ whose mass values are non-zero are called \emph{focal elements} of $m$.

The \emph{belief function} (BF) associated with a BPA $m : 2^\Theta\rightarrow[0,1]$ is the set function 
\[
Bel(A) = \sum_{B\subseteq A} m(B). 
\]
The corresponding \emph{plausibility function} is 
\[
Pl(A) \doteq \sum_{B\cap A\neq \emptyset} m(B) \geq Bel(A).
\]
Belief functions can be defined axiomatically, however, without resorting to the above random set interpretation \cite{Shafer76}. 

Classical probability measures on $\Theta$ are a special case of belief functions (those assigning mass to singletons only), termed \emph{Bayesian belief functions}.
A BF is said to be \emph{consonant} if its focal elements $A_1,...,A_m$ are nested: $A_1 \subset \cdots \subset A_m$, and corresponds to a possibility measure.

\subsubsection{Combination}

\emph{Dempster's combination} $Bel_1 \oplus Bel_2$ 
of two belief functions 
on $\Theta$ is the unique BF there with as focal elements all the {non-empty} intersections of focal elements of $Bel_1$ and $Bel_2$, and basic probability assignment
\[
m_{\oplus}(A) = \frac{m_\cap(A)} {1- m_\cap(\emptyset)},
\]
where 
\[
m_\cap(A) = \sum_{B \cap C = A} m_1(B) m_2(C) 
\]
and $m_i$ is the BPA of the input BF $Bel_i$.

Dempster's combination naturally induces a conditioning operator. Given a conditioning event $A \subset \Theta$, the `logical' or `categorical' belief function $Bel_A$ such that $m(A)=1$ is combined via Dempster's rule with the a-priori belief function $Bel$. The resulting BF $Bel \oplus Bel_A$ is the {conditional belief function given $A$} \emph{a la Dempster}, denoted by $Bel_\oplus(A|B)$.

Many alternative combination rules have since been defined \cite{Klawonn:1992:DBT:2074540.2074558,yager87on,dubois88representation,DENOEUX2008234}, often associated with a distinct approach to conditioning \cite{Denneberg1994,fagin91new,suppes1977}. An exhaustive review of these proposals can be found in \cite{cuzzolin2021springer}, Section 4.3. 

\subsubsection{Belief functions and other measures}

Each belief function $Bel$ uniquely identifies a credal set \cite{kyburg87bayesian}
\[
\mathcal{P}[Bel] = \{ P \in \mathcal{P} : P(A) \geq Bel(A) \}
\]
(where $\mathcal{P}$ is the set of all probabilities one can define on $\Theta$), of which it is its lower envelope:
$Bel(A) = \underline{P}(A)$. 
The probability intervals resulting from Dempster's updating of the credal set associated with a BF, however, are included in those resulting from Bayesian updating \cite{kyburg87bayesian}.

Belief functions are infinitely monotone capacities, and a special case of coherent lower previsions. 
On the other hand, possibility measures are equivalent to consonant BFs \cite{cuzzolin2021springer}. Belief functions can also be generalised to assume values on fuzzy sets, rather than traditional `crisp' ones \cite{Biacino07}.
\\
Indeed, when we recall that belief measures are induced by mass assignments on subsets of a frame we realise that {generalisations of belief functions which are defined on fuzzy sets} can (and indeed have) been proposed .
Following Zadeh's work, Ishizuka et al \cite{ISHIZUKA1982179}, Ogawa and Fu \cite{OGAWA1985295}, Yager \cite{YAGER198245}, Yen \cite{yen90generalizing} and recently Biacino \cite{Biacino07} extended belief theory to fuzzy sets 
by defining an appropriate measure of inclusion for them.
Every BF specifies a unique p-box, whereas each p-box specifies an entire equivalence class of belief functions \cite{REGAN20041}.

Authors such as Heilpern \cite{heilpern97representation}, Yager \cite{yager99-class,INT:INT4550010106}, Palacharla \cite{palacharla94understanding}, Romer \cite{romer95applicability}, Kreinovich \cite{Kreinovich01tech}, Denoeux \cite{denoeuxmodeling} and many others \cite{Renaud99,goodman82fuzzy,Feng201287,781806,Florea03,sevastianov07-numerical} have also studied the connection between fuzzy and Dempster--Shafer theory. 

A framework similar to belief functions' is Spohn's, 
in which some propositions are believed to be true, others are believed to be false, and the remainder are neither believed nor disbelieved \cite{Spohn1988}.

\subsubsection{Behavioural interpretation}

An intriguing behavioural interpretation of belief functions has been recently advanced in \cite{Kerkvliet2017}. A betting function $R : \mathcal{L} \rightarrow [0,1]$ is a binary function on the set of gambles such that $\forall X \in \mathcal{L}$ $\exists \alpha_X \in \mathbb{R}$ such that $R(X + \alpha) = 0$ for $\alpha < \alpha_X$, while $R(X + \alpha) =1$ for $\alpha \geq \alpha_X$. 

Given a betting function, 
\[
Buy_R(X) \doteq \max \{ \alpha \in \mathbb{R} : R(X - \alpha) =1 \} 
\]
is the maximum price an agent is willing to pay for the gamble $X$. In Walley's terminology, the betting function determines the set of desirable gambles $\mathcal{D} = \{X : R(X) = 1\}$, whereas $Buy_R(X)$ coincides with  (the lower prevision $\underline{P}(X)$ of $X$). 

A \emph{belief valuation} $\mathcal{B}$ is a belief function $Bel_B$, such that $Bel_B(A) = 1$ if $A \supseteq B$, 0 otherwise.
For any belief valuation $\mathcal{B}$, the \emph{guaranteed revenue} $G_\mathcal{B} : \mathcal{L} \rightarrow \mathbb{R}$ is defined as
\[
G_\mathcal{B}(X) \doteq \max_{A \supseteq B} \min_{\theta \in A} X(\theta).
\]
A betting function $R$ is \emph{B-consistent} if for all $X_1,...,X_N \in \mathcal{L}$, $Y_1,...,Y_M \in \mathcal{L}$ such that 
\[
\sum_{i=1}^N G_\mathcal{B} (X_i) \leq \sum_{j=1}^M G_\mathcal{B} (X_j) 
\]
for every belief valuation $\mathcal{B}$, we have 
\[
\sum_{i=1}^N Buy_R (X_i) \leq \sum_{j=1}^M Buy_R (X_j).
\] 

$Bel$ is a belief function if and only if there exists a coherent (in Walley's sense)
and B-consistent $R$ such that $Bel(A) = Buy_R(1_A)$ \cite{Kerkvliet2017}. 
Also, if $R$ is a coherent betting function, then $R$ is B-consistent if and only if there is a BPA $m$ such that: $Buy_R(X) = \sum_{A \subseteq \Theta} m(A) \min_{\theta \in A} X(\theta)$ for all $X \in \mathcal{L}$.
In other words, adding B-consistency to imprecise probabilities' rationality axioms, lower previsions reduce to belief functions.

\subsection{Random sets} \label{sec:random-sets}

When the sample space is continuous, belief functions generalise to \emph{random sets}.

Let $(\Omega, \mathcal{F},P)$ be a probability space \cite{Molchanov05}.
A map $X : \Omega \rightarrow \mathcal{C}$, where $\mathcal{C}$ is the space of closed subsets of a certain topological space $\mathbb{E}$, is called a \emph{random closed set} if, for every compact set $K \in \mathcal{K}$ in $\mathbb{E}$
\[
\{ \omega : X(\omega) \cap K \neq \emptyset \} \in \mathcal{F}.
\]
Alternatively, we can define the Borel $\sigma$-algebra $ \mathcal{B}(\mathcal{C})$ generated by $\{C \in \mathcal{C} : C \cap X \neq \emptyset\}$ for all compact sets $K$, and call
$X$ a random closed set if 
it is measurable with respect to the Borel $\sigma$-algebra on $\mathcal{C}$ with respect to the Fell topology, namely:
\[
X^{-1}(\mathcal{X}) = \{ \omega : X(\omega) \in \mathcal{X} \} \in \mathcal{F}
\]
for each $\mathcal{X} \in \mathcal{B}(\mathcal{C})$.
The distribution of a random closed set $X$ is determined by $P(\mathcal{X}) = P(\{\omega : X(\omega) \in \mathcal{X}\})$ for all $\mathcal{X} \in \mathcal{B}(\mathcal{C})$. 
In particular we can consider the case in which $\mathcal{X} = \mathcal{C}_K \doteq \{C \in \mathcal{C} : C \cap K \neq \emptyset \}$ and $ P(\{X(\omega) \in \mathcal{C}_K \}) = P(\{X \cap K \neq \emptyset \})$, since the families $\mathcal{C}_K, K\in\mathcal{K}$, generate the Borel $\sigma$-algebra $\mathcal{B}(\mathcal{C})$.

A functional $T_X : \mathcal{K} \rightarrow [0, 1]$ given by
\[
T_X (K) = P (\{X \cap K \neq \emptyset \}), 
\quad
K \in \mathcal{K}
\]
is termed the \emph{capacity functional} of $X$. 
In particular, if $X = \{\xi\}$ is a random singleton, i.e., a classical random variable, then $T_X(K) = P(\{\xi \in K\})$, so that the capacity functional is the probability distribution of the RV $\xi$.

The functional $T_X$ is indeed a \emph{capacity} in the sense of Section \ref{sec:others-capacities}.
The links betwen random closed sets and belief functions, upper and lower probabilities, and contaminated models in statistics is briefly hinted at in \cite{Molchanov05}, Chapter 1.

\section{Measures on functional spaces} \label{sec:others-hop}

Possibly the most general group of approaches clusters those advocating the need to go beyond the notion of event entirely. These methods replace the usual $\sigma$-algebra of events with a functional space, the space of all functions of a certain class that can be defined over the universe of discourse $\Theta$, and attach numerical quantifiers to these functions. Lower and upper previsions (see Section \ref{sec:walley}) are an example of such quantifiers, as they are defined on gambles rather than events.

\subsection{Higher-order probabilities}

\subsubsection{Second-order probabilities}

An intuitive approach is to model uncertainty about probabilities using probability itself \cite{Domotor81}.
As pointed out by \cite{Baron1987}, a \emph{second-order probability} $Q(P)$ may be understood as the probability that the `true' probability of something has the value $P$. 
Baron proceeded to derive a rule for combining evidence from two independent sources generating each a second-order probability $Q_i(P)$, and showed that Dempster's rule is a special case representing a restriction of a full Bayesian analysis.
In Fung and Chong's \emph{metaprobability theory}, belief update is still based on Bayes' rule \cite{DBLP:journals/corr/abs-1304-3427}, namely 
\[
p^2 (p| D, Pr) \propto p^2 (D | p, Pr) \cdot p^2(p | Pr), 
\]
where $D$ is the data (evidence) and $Pr$ a prior on space of (first order) probability distributions $p$.

\subsubsection{Josang's Dirichlet distributions}

Josang showed that a bijective mapping exists between Dirichlet distributions \cite{Josang2007}, a special class of second-order probability distributions, and belief functions whose focal elements are of size 1 or $|\Theta|$. 
Such a link can be exploited to apply belief based reasoning to statistical data, or to apply statistical and probabilistic analysis to belief functions.

\subsubsection{Gaifman's higher order probability spaces}

\cite{Gaifman1988} axiomatically defined a \emph{higher order probability space} (HOP) as consisting of a probability space and an operation $PR$, such that, for every event $A$ and every real closed interval $\Delta$, $PR(A ,\Delta)$ is the event that $A$'s `true' probability (in the sense described above) lies in $\Delta$. 
In a \emph{general} HOP the operation $PR$ includes an additional argument ranging over an ordered set of time-points, so that $PR(A,t,\Delta)$ is the event that $A$'s probability at time $t$ lies in $\Delta$. Various connections with modal logic were pointed out.

\subsubsection{Kyburg's analysis}

Whereas
\cite{Skyrms1980-SKYHOD} argued that higher order probabilities are essential for a correct representation of belief,
\cite{kyburg88hop} claimed that higher order probabilities can always be replaced by the
marginal distributions of a joint probability distribution defined on $I \times \Theta$, where $\Theta$ is the sample space of lower order probabilities, and $I$ is a finite set indexing all possible probability distributions $P_i$ (i.e., the sample space for the second-order probability $Q$).

This follows from the principle that 
\cite{Skyrms1980-SKYHOD}
\[
P(\theta) = \sum_{i \in I} Q(P_i) P_i(\theta) = E[P_i(\theta)],
\]
combined with the use of expected utilities.
As a result, Kyburg claimed, higher order probabilities do not appear to offer any conceptual or computational advantage.

\begin{figure*}[ht!]
\begin{center}
\includegraphics[width=0.85\textwidth]{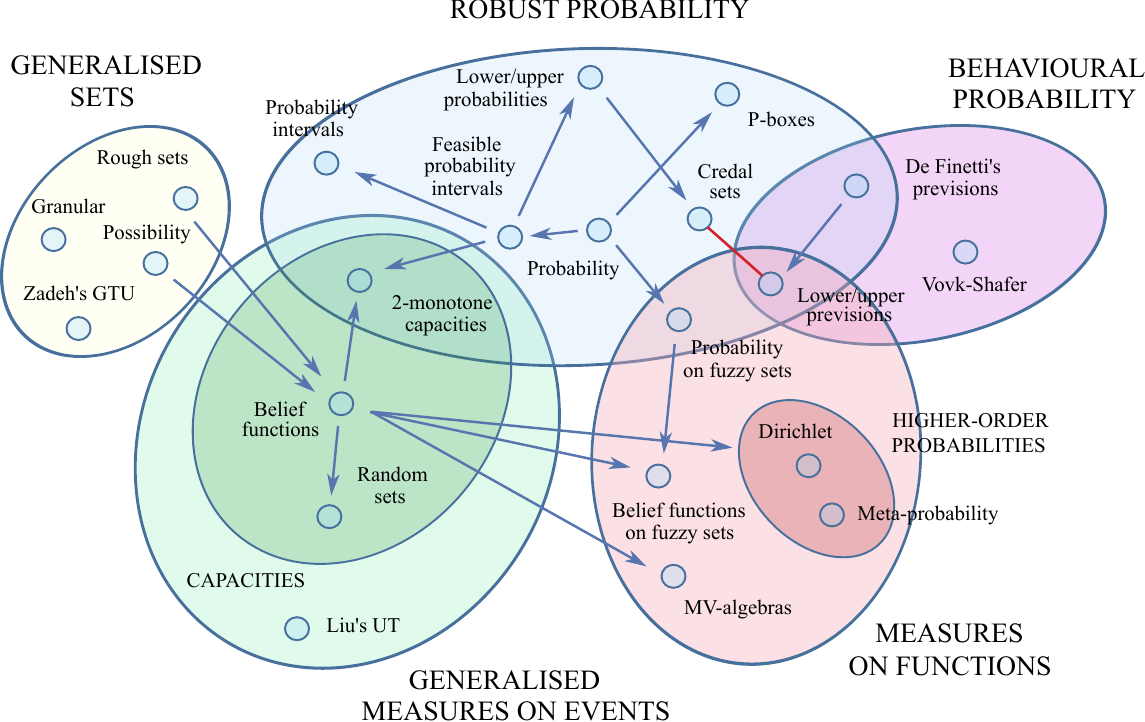}
\end{center}
\caption{Illustration of how the uncertainty theories reviewed in this paper cluster into a number of groups, associated with the rationale and class of objects they attempt to mathematically model. Arrows represent implications, from the less general model to the more general. 
 \label{fig:diagram}} 
\end{figure*}

\subsection{MV-algebras} 

Various extensions of probability measures of belief functions to \emph{fuzzy events}, rather than traditional, `crisp' events, have also been proposed. As fuzzy events amount to membership functions on the universe $\Omega$, they can also be assimilated to measures on functional spaces.
Another remarkable example is provided by {MV-algebras} \cite{Kroupa2010,Flaminio2011}. 

An \emph{MV-algebra} is an algebra $\langle M, \oplus, \neg, 0 \rangle$ of many-valued events, upon which upper and lower probabilities can be defined. $M$ is endowed with a binary operation $\oplus$, a unary operation $\neg$ and a constant 0 such that $\langle M, \oplus, 0 \rangle$ is an Abelian monoid and the following equations hold true for every $f ,g \in M$:
$\neg \neg f = f$, $f \oplus \neg 0 = \neg 0$, $\neg (\neg f \oplus g) \oplus g = \neg (\neg g \oplus f) \oplus f$. Boolean algebras are, in fact, a special case of MV algebras. 

A \emph{state} is a mapping $s:M \rightarrow [0,1]$ such that $s(1)=1$ and 
\[
s(f+g) = s(f)+s(g) \; \text{whenever} \; f \odot g = 0, 
\]
where $f \odot g \doteq \neg (\neg f \oplus \neg g)$. Clearly, states are generalisations of finitely additive probability measures, once we replace events $A,B$ with continuous functions $f,g$ onto $[0,1]$ and $\cap$ with $\odot$. 

Considering the MV algebra $[0,1]^{2^\Theta}$ of all functions $2^\Theta \rightarrow [0,1]$, where $X$ is finite, one define a belief function $Bel : [0,1]^{2^\Theta} \rightarrow [0,1]$ there if there is a state on the MV algebra $[0,1]^{\mathcal{P}(X)}$ such that $s(1_{\emptyset}) = 0$ and $Bel(f) = s(\rho(f))$, for every $f \in [0,1]^X$. 

\section{Discussion}

\subsection{A tentative classification} 

Figure \ref{fig:diagram} provides a visual illustration of our tentative classification of uncertainty theories into clusters of approaches, depending on what objects they attach values to and, related to this, their rationale. The diagram also shows the dependency between the various theories, based on the results available so far. An arrow between formalism 1 and formalism 2 means that the former is less general than the latter. The `quality' of the various proposals differs radically.
For instance, whereas an uncertain entropy and an uncertain calculus are built within Liu's UT, the general lack of rigour and convincing justification leaves us quite unimpressed with this work.
As far as Zadeh's GTU is concerned, which was not reviewed here in detail, generality is achieved there in a rather nomenclative way, which explains the complexity and lack of clear rationale for the formalism.

Imprecise probability is remarkable for its generality, as special cases of coherent lower previsions include probabilities, de Finetti previsions, Choquet capacities, possibility and belief measures, random sets but also probability boxes, (lower and upper envelopes of) credal sets, and robust Bayesian models. Nevertheless, IP has so far failed to achieved the expected impact in both mathematical statistics and application fields. This may lead to the reflection that generality, after all, is not the real thing. As discussed, whether very general settings such as meta-probability offer any conceptual or practical advantage is unclear. It is also somewhat puzzling that coherent lower previsions end up being as expressive as credal sets, which can be handled without resorting to complex linear programming optimisation in a robust Bayesian fashion. Moreover, the linear utility assumptions in IP are very strong, which weakens their case as a `natural' formalisation of subjective reasoning.

Other recent and philosophically interesting frameworks such as Vovk and Shafer's game theoretical probability (also related to imprecise probability, \cite{miranda2008survey}) have encountered similar issues with rationale and motivation. In addition, to this author's understanding, a further extension of this reasoning to non-repetitive situations is still sought.

\subsection{Geometry as unifying language}

\emph{Geometry} may be a unifying language for the field \cite{cuzzolin18belief-maxent,Cuzzolin99,cuzzolin05isipta,cuzzolin00mtns,cuzzolin13fusion,gennari02-integrating,gong2017belief,black97geometric,rota97book,ha98geometric,wang91geometrical}, possibly in conjunction with an algebraic view \cite{cuzzolin00rss,cuzzolin01bcc,cuzzolin08isaim-matroid,cuzzolin01lattice,cuzzolin05amai,cuzzolin07bcc,cuzzolin14algebraic}.

In the {geometric approach to uncertainty}, uncertainty measures can be seen as points of a suitably complex geometric space, and there manipulated (e.g. combined, conditioned and so on) \cite{cuzzolin01thesis,cuzzolin2008geometric,cuzzolin2021springer}.
Much work has been focusing on the geometry of belief functions, which live in a convex space termed the \emph{belief space}, which can be described both in terms of a simplex (a higher-dimensional triangle) and in terms of a recursive bundle structure \cite{cuzzolin01space,cuzzolin03isipta,cuzzolin14annals,cuzzolin14lap}. The analysis can be extended to Dempster's rule of combination by introducing the notion of a conditional subspace and outlining a geometric construction for Dempster's sum \cite{cuzzolin02fsdk,cuzzolin04smcb}.
The combinatorial properties of plausibility and commonality functions, as equivalent representations of the evidence carried by a belief function, have also been studied \cite{cuzzolin08pricai-moebius,cuzzolin10ida}. The corresponding spaces are simplices which are congruent to the belief space.
\\
Subsequent work extended the geometric approach to other uncertainty measures, focusing in particular on possibility measures (consonant belief functions) \cite{cuzzolin10fss} and consistent belief functions \cite{cuzzolin11-consistent,cuzzolin09isipta-consistent,cuzzolin08isaim-simplicial}, in terms of simplicial complexes \cite{cuzzolin04ipmu}. Analyses of belief functions in terms credal sets have also been conducted \cite{cuzzolin08-credal,antonucci10-credal,burger10brest}.

The geometry of the relationship between measures of different kinds has also been extensively studied \cite{cuzzolin05hawaii,cuzzolin09-intersection,cuzzolin07ecsqaru,cuzzolin2010credal}, with particular attention to the problem of transforming a belief function into a classical probability measure \cite{Cobb03isf,voorbraak89efficient,Smets:1990:CPP:647232.719592}. One can distinguish between an `affine' family of probability transformations \cite{cuzzolin07smcb} (those which commute with affine combination in the belief space), and an `epistemic' family of transforms \cite{cuzzolin07report}, formed by the relative belief and relative plausibility of singletons \cite{cuzzolin08unclog-semantics,cuzzolin2008semantics,CUZZOLIN2012786,cuzzolin06-geometry,cuzzolin10amai}, which possess dual properties with respect to Dempster's sum \cite{cuzzolin2008dual}.
The problem of finding the possibility measure which best approximates a given belief function \cite{aregui08constructing} can also be approached in geometric terms \cite{cuzzolin09ecsqaru,cuzzolin11isipta-consonant,cuzzolin14lp,Cuzzolin2014tfs}. In particular, approximations induced by classical Minkowski norms can be derived and compared with classical outer consonant approximations \cite{Dubois90}.
Minkowski consistent approximations of belief functions in both the mass and the belief space representations can also be derived \cite{cuzzolin11-consistent}.

In fact, the geometric approach to uncertainty can be applied to various other elements of the inference chain. For instance, the conditioning problem can be posed geometrically \cite{lehrer05updating}. Conditional belief functions can be defined as those which minimise an appropriate distance between the original belief function and the `conditioning simplex' associated with the conditioning event \cite{cuzzolin10brest,cuzzolin11isipta-conditional}. 
A semantics for the main probability transforms can be provided in terms of credal sets, i.e., convex sets of probabilities \cite{cuzzolin2010credal}. 

Recent papers on this topic include \cite{luo2020vector,pan2020probability,long2021visualization}.

\subsection{What is the `right' formalism to use?} 

A current of thought supports the notion
that there is no such a thing as \emph{the} best mathematical model \cite{Fagin88,klir95principles}. The choice of the most suitable methodology, instead, should depend on the actual problem at hand. 
We are not quite in favour of such a position. In our view, evidence from the applications should inform the choice of the most suitable formalism, general enough to provide a comprehensive theory of uncertainty, but no so general as to become useless or overcomplicated. Uncertainty theory should avoid resorting to axioms whenever possible, as this inevitably opens the way to arguments on their validity, especially when claims are made regarding human rationality.

In this sense, strong evidence that observations are inherently set-valued provides support for the theory of random sets. The latter is very general, belongs to most of the above conceptual clusters (amounts to a robust probability theory, generalises set and measure theory as well as Bayesian reasoning, has significant behavioural interpretations), and naturally encodes ignorance, missing, vague and propositional data, and is compatible with both Bayesian and frequentist interpretations of probability \cite{cuzzolin2021springer}. Theoretical advances have been made in recent times, but more are needed.

Despite its flaws additive probability is simple to understand and easy to use. To overcome the conservatorism of experts and practitioners alike, uncertainty theory needs to provide a similar compelling, simple message.

\bibliographystyle{plain}
\bibliography{arxiv-ijcai-survey}

\end{document}